\documentclass{article}
\usepackage{amssymb}
\title{\textsf{The Finitary Andrews-Curtis Conjecture}}

\author{\textsf{ Alexandre V. Borovik} \and \textsf{Alexander Lubotzky} \and
\textsf{ Alexei G. Myasnikov}}
\date{\textsf{12 June 2003}}

\newtheorem{theorem}{Theorem}[section]
\newtheorem{lemma}[theorem]{Lemma}
\newtheorem{fact}[theorem]{Fact}
\newtheorem{corollary}[theorem]{Corollary}

\newtheorem{definition}{Definition}

\newcommand{\bi}{\begin{itemize}}
\newcommand{\ei}{\end{itemize}}
\newcommand{\rmit}[1]{\item[{\rm #1}]}

\begin{document}

\pagestyle{myheadings} \markright{{\scriptsize \textsf{A.V.
Borovik, A. Lubotzky and A.G. Myasnikov} $\bullet$ \textsf{ The
Finitary Andrews-Curtis Conjecture} $\bullet$ 12.06.03}}

\maketitle

\begin{center}
{\large\emph{To Slava Grigorchuk as a token of our friendship.}}
\end{center}

\bigskip

\begin{abstract}
The well known Andrews-Curtis Conjecture \cite{AC} is still open.
In this paper, we establish its finite version by describing
precisely the connected components of the Andrews-Curtis graphs
of finite groups. This finite version has independent importance
for computational group theory. It also resolves a question asked
in \cite{BKM} and shows that a computation in finite groups
cannot lead to a counterexample to the classical conjecture, as
suggested in \cite{BKM}.

\end{abstract}

\section{Andrews-Curtis graphs}
\label{se:2-1}

Let $G$ be a group and  $G^k$ be the set of all $k$-tuples of
elements of $G$.

The following  transformations of the set $G^k$ are called {\it
elementary Nielsen transformations {\rm (}or moves{\rm )}}:
\begin{enumerate}
\item [(1)] $(x_1, \ldots, x_i, \ldots,  x_k) \longrightarrow (x_1,
\ldots,x_ix_j^{\pm 1}, \ldots,  x_k), \; i \ne j;$
 \item [(2)] $(x_1, \ldots, x_i, \ldots, x_k)  \longrightarrow  (x_1,
\ldots,x_j^{\pm 1}x_i, \ldots,  x_k), \; i \ne j;$
 \item [(3)] $(x_1, \ldots, x_i, \ldots, x_k)  \longrightarrow  (x_1,
\ldots,x_i^{-1}, \ldots, x_k).$
\end{enumerate}

Elementary Nielsen moves transform generating tuples of $G$ into
generating tuples.  These moves together with the transformations
\begin{enumerate}
 \item [(4)] $(x_1, \ldots, x_i, \ldots, x_k) \longrightarrow
(x_1,\ldots,x_i^{w}, \ldots, x_k ), \;  w\in S \cup S^{-1} \subset
G,$
\end{enumerate}
where $S$ is a fixed subset of $G$, form a set of {\it elementary
Andrews-Curtis transformations relative to $S$} (or, shortly,
$AC_S$-moves). If $S=G$ then AC-moves transform n-generating
tuples (i.e., tuples which generate $G$ as a normal subgroup) into
n-generating tuples. We say that two $k$-tuples $U$ and $V$ are
$AC_S$-equivalent, and write $U \sim_S V$, if there is a finite
sequence of $AC_S$-moves which transforms $U$ into $V$. Clearly,
$\sim_S$ is an equivalence relation on the set $G^k$ of $k$-tuples
of elements from $G$. In the case when $S = G$ we omit $S$ in the
notations and refer to $AC_S$-moves simply as to $AC$-moves.

 We slightly change notation from that of
\cite{BKM}. For a subset $Y \subset G$ we denote by $gp_{G}(Y)$
the normal closure of $Y$ in $G$,  by $d(G)$   the minimal number
of generators of $G$, and by $d_G(G)$  the minimal number of
normal generators of $G$. Now, $d_G(G)$ coincides with $nd(G)$ of
\cite{BKM}.

Let  $N_k(G)$, $k \geqslant d_G(G),$ be the set of all $k$-tuples
of elements in $G$ which generate $G$ as a normal subgroup:
$$ N_k(G) =  \left\{\,
(g_1,\ldots, g_k) \mid gp_G(g_1,\ldots, g_k) = G \,\right\}.$$
Then the {\em Andrews--Curtis graph\/} $\Delta_k^S(G)$ of the
group $G$  with respect to a given subset $S \subset G$ is the
graph whose vertices are $k$-tuples from $N_k(G)$ and such that
two vertices are connected by an edge if one of them is obtained
from another by an elementary $AC_S$-transformation.  Again, if $S
= G$ then we refer to $\Delta_k^G(G)$ as to the Andrews-Curtis
graph of $G$ and denote it by $\Delta_k(G)$. Clearly, if $S$ is a
generating set of $G$ then  the graph $\Delta_k^S(G)$ is connected
if and only if the graph $\Delta_k(G)$ is connected. Observe, that
if $S$ is finite then $\Delta_k^S(G)$ is a regular graph of finite
degree.

 The famous Andrews-Curtis
conjecture \cite{AC} can be stated in the following way.

\begin{quote} \textsf{\sl {\bf AC-Conjecture:}   For a
free group $F_k$ of rank \/ $k \geqslant 2$, the Andrews--Curtis
graph\/ $\Delta_k(F_k)$ is connected.}
\end{quote}

There are some doubts whether this well known old conjecture is
true. Indeed, Akbulut and Kirby \cite{AK} suggested a series of
potential counterexamples for $k = 2$:
 \begin{equation}
 \label{eq:AC}
 (u, v_n) = (xyxy^{-1}x^{-1}y^{-1}, x^ny^{-(n+1)}), \ \ n\geqslant 2.
 \end{equation}

  In \cite{BKM}, it has been suggested that one may be able to confirm
one of these potential counterexamples by showing that for some
homomorphism $\phi:F_2 \rightarrow G$ into a  finite group $G$ the
pairs  $(u^\phi,v_n^\phi)$ and $(x^\phi,y^\phi)$  lie in different
connected components of $\Delta_2(G)$. Notice that in view of
\cite{My}\label{My-mentioned} the group $G$ in the counterexample
cannot be soluble.

Our main result describes the connected components of the
Andrews-Curtis graph of a finite group. As a corollary we show
that $(u^\phi,v_n^\phi)$ and $(x^\phi,y^\phi)$  lie in the same
connected components of $\Delta_2(G)$ for every finite group $G$
and any homomorphism $\phi:F_2 \rightarrow G$, thus resolving the
question from \cite{BKM}.

\begin{theorem}
Let $G$ be a finite group and $k \geqslant \max\{d_G(G),2\}$. Then
two tuples $U, V$ from $N_k(G)$ are AC-equivalent if and only if
they are AC-equivalent in the abelianisation ${\rm Ab}(G) =
G/[G,G]$, i.e., the connected components of the AC-graph
$\Delta_k(G)$ are precisely the preimages of the connected
components of the AC-graph $\Delta_k({\rm Ab}(G))$.
\label{th:G=Omega} \label{th:lifting}
\end{theorem}

Notice that, for the abelian group $A = {\rm Ab}(G)$, a normal
generating set is just a generating set and the non-trivial
Andrews-Curtis transformations are Nielsen moves (1)--(3).
Therefore the vertices of $\Delta_k({A})$ are the same as these of
the \emph{product replacement graph} $\Gamma_k(A)$ \cite{dg,P2}:
they are all generating $k$-tuples of $A$. The only difference
between $\Gamma_k(A)$ and $\Delta_k(A)$ is that the former has
edges defined only by `transvections' (1)--(2), while in the
latter the inversion of components (3) is also allowed. The
connected components of product replacements graphs $\Gamma_k(A)$
for finite abelian groups $A$ have been described by Diaconis and
Graham \cite{dg}; a slight modification of their proof leads to
the following observation

\begin{fact}[Diaconis and Graham \cite{dg}] Let\/ $A$ be a finite abelian
group and
$$
A = Z_1 \times\cdots\times Z_d
$$
its canonical decomposition into a direct product of cyclic groups
such that\/ $|Z_i|$ divides $|Z_j|$ for\/ $i < j$. Then \bi
\rmit{(a)} If\/ $k > d$ then $\Delta_k(A)$ is connected.
\rmit{(b)} If\/ $k=d \geqslant 2$, fix generators $z_1,\dots,z_d$
of the subgroups $Z_1,\dots,Z_d$, correspondingly. Let\/ $m =
|Z_1|$. Then $\Delta_d(A)$ has $\phi(m)/2$ connected components
{\rm (}here $\phi(n)$ is the Euler function{\rm )}. Each of these
components has a representative of the form
$$
(z_1^\lambda, z_2,\dots, z_d), \quad \lambda \in
({\mathbb{Z}}/m{\mathbb{Z})^*}.
$$
Two tuples
$$(z_1^\lambda, z_2,\dots, z_d) \;\hbox{ and }\; (z_1^\mu,
z_2,\dots, z_d),\quad \lambda, \mu \in
({\mathbb{Z}}/m{\mathbb{Z}})^*,$$ belong to the same connected
component if and only if\/ $\lambda = \pm \mu$. \ei
\label{fact:dg}
\end{fact}

Taken together, Theorem~\ref{th:lifting} and Fact~\ref{fact:dg}
give a complete description of components  of the Andrews-Curttis
graph $\Delta_k({G})$ of a finite group $G$.

Notice that in an abelian group $A$ \begin{eqnarray*}
(xyxy^{-1}x^{-1}y^{-1}, x^ny^{-(n+1)}) &\sim & (xy^{-1},
x^ny^{-(n+1)})\\
& \sim & (xy^{-1}, x^{n-1}y^{-n})\\
& \vdots &\\
& \sim & (yx^{-1},y^{-1})\\
& \sim & (x,y)
\end{eqnarray*}
 so for every homomorphism $\phi:F_2 \rightarrow G$ as above the
 images $(u^\phi,v_n^\phi)$ and $(x^\phi,y^\phi)$ are AC equivalent in the
 abelianisation  of $G$, hence they lie in the
 same connected component of $\Delta_2(G)$.

The following corollary of Theorem~\ref{th:lifting} leaves no
hope of finding an counterexample to the Andrews-Curtis
conjecture by looking at the connected components of the
Andrews-Curtis graphs of finite groups.

\begin{corollary}
For any $k \geqslant 2$, and any epimorphism $\phi: F_k
\rightarrow G$ onto a finite group $G$, the image of\/
$\Delta_k(F_k)$ in $\Delta_k(G)$ is connected.
\end{corollary}

One may try to reject the AC-conjecture by testing AC-equivalence
of the tuples $(u,v_n)$ and $(x,y)$ in the {\em infinite}
quotients of the group $F_2$. To this end we introduce the
following definition.

\begin{quote}
 \textsf{\sl {\bf Definition:} \/ We say that a group $G$ satisfies the {\em
generalised
Andrews-Curtis conjecture} if for any $k \geqslant
\max\{\,d_G(G),2\,\}$ tuples $U, V \in N_k(G)$ are AC-equivalent
in $G$ if and only if their images are AC-equivalent in the
abelianisation ${\rm Ab}(G)$.}
\end{quote}

\begin{quote}
 \textsf{\sl {\bf Problem:} \/ Find a group $G$ which does not satisfy the
 generalised Andrews-Curtis conjecture.}
\end{quote}

It will be interesting to look, for example, at the Grigorchuk
group \cite{rg1,rg2}. It is a finitely generated residually finite
$2$-group $G$ which is just-infinite, that is, every normal
subgroup has finite index. Therefore the generalised
Andrews-Curtis conjecture holds in every proper factor group of
$G$ by Theorem~\ref{th:lifting}. What might be also relevant, the
conjugacy problem in the Grigorchuk group is solvable
\cite{leonov,rozhkov,BGS}. This makes the Grigorchuk group a very
interesting testing ground for the generalised Andrews-Curtis
conjecture.

\section{Relativised  Andrews-Curtis graphs and black-box groups}
\label{se:relAC}

Following \cite{BKM}, we also introduce a relativised version of
the Andrews-Curtis transformations of the set $G^k$ for the
situation when $G$ admits some fixed group of operators  $\Omega$
(that is, a group $\Omega$ which acts on $G$ by automorphisms); we
shall say in this situation that $G$ is an {\em
$\Omega$-group}\footnote{We shall use the terms $\Omega$-subgroup,
normal $\Omega$-subgroup, $\Omega$-simple $\Omega$-subgroup, etc.
in their obvious meaning.}.  In that case, we view the group $G$
as a subgroup of the natural semidirect product $G\cdot \Omega$ of
$G$ and $\Omega$. In particular, the set of $AC_{G \Omega}$-moves
is defined and the set $G^k$ is invariant  under these moves.  In
particular, if $N$ is a normal subgroup of $G$, we view $N$ as a
$G$-subgroup in the sense of this definition. As we shall soon
see, $AC_{G\Omega}$-moves appear in the product replacement
algorithm for generating pseudo-random elements of a normal
subgroup in a black box finite group.

For a subset  $Y \subset G$ of an $\Omega$-group $G$ we denote
 by $gp_{G\Omega}(Y)$  the normal
closure of $Y$ in $G\cdot \Omega$,  and by $d_{G\Omega}(G)$ the
minimal number of normal generators of $G$ as a normal subgroup of
$G\cdot \Omega$.

Let  $N_k(G,\Omega)$, $k \geqslant d_{G\Omega}(G),$ be the set of
all $k$-tuples of elements in $G$ which generate $G$ as a normal
$\Omega$-subgroup:
$$ N_k(G,\Omega) =  \left\{\,
(g_1,\ldots, g_k) \mid gp_{G\cdot \Omega}(g_1,\ldots, g_k) = G
\,\right\}.$$
 Then the {\em relativised Andrews--Curtis graph\/}
$\Delta_k^\Omega(G)$ of the group $G$  is the graph whose vertices
are $k$-tuples from $N_k(G,\Omega)$ and such that two vertices are
connected by an edge if one of them is obtained from another by an
elementary $AC_{G\Omega}$-transformation.

A {\em black box group} $G$ is a finite group with a device (`oracle')
 which produces
its (pseudo)random (almost) uniformly distributed elements; this
concept is of crucial importance for computational group theory,
see  \cite{kantor}. If the group $G$ is given by generators, the
so-called {\em product replacement algorithm} \cite{celler2,P2}
provides a very efficient and practical way of producing random
elements from $G$; see \cite{LuP} for a likely theoretical
explanation of this (still largely empirical) phenomenon in terms
of the (conjectural) Kazhdan's property (T) \cite{K} for the group
of automorphisms of the free group $F_k$ for $k >4$. In the
important case of generation of random elements in a normal
subgroup $G$ of a black box group $\Omega$, the following simple
procedure is a modification of the product replacement algorithm:
start with the given tuple $U \in N_k(G,\Omega)$, walk randomly
over the graph $\Delta_k^\Omega(G)$ (using the `oracle' for
$\Omega$ for generating random $AC_{G\Omega}$-moves   and return
randomly chosen components $v_i$ of vertices $V$ on your way. See
\cite{borovik,BKM,L-GM} for a more detailed discussion of this
algorithm, as well as its further enhancements.

Therefore the understanding of the structure---and ergodic properties---of
the Andrews-Curtis graphs
$\Delta_k^\Omega(G)$ is of some importance for the theory of black box
groups.

\medskip
The following results are concerned with the connectivity of the
relativised  Andrews-Curtis graphs of finite groups.

\begin{theorem} Let $G$ be a finite $\Omega$-group
which is perfect as an abstract group, $G = [G,G]$.
Then the graph $\Delta_k^\Omega(G)$ is connected for every
$k\geqslant 2$.
\label{normal}
\end{theorem}

Of course, this result can be immediately reformulated for normal
subgroups of finite groups:

\begin{corollary} Let $G$ be a finite group
and  $N \triangleleft G$ a perfect normal subgroup. Then the graph
$\Delta_k^G(N)$ is connected for every $k\geqslant 2$.
\label{cor:normal}
\end{corollary}

 We would like to record another immediate corollary of
Theorem~\ref{normal}.

\begin{corollary}
Let $G$ be a perfect finite group, $g_1,\ldots, g_k$, $k \geqslant 2$
generate
$G$ as a normal subgroup and $\phi: F_k \longrightarrow G$ an epimorphism.
Then there exist $f_1,\ldots,f_k \in F_k$ such that $\phi(f_i) = g_i$, $i
=1,\ldots, k$,
and $f_1,\ldots, f_k$ generate $F_k$ as a normal subgroup.
\end{corollary}

Note that if we take $g_1,\dots,g_k$ as a set of generators for
$G$, then in general we cannot pull them back to a set
$f_1,\dots, f_k$ of generators for $F_k$, an example can be found
in $G = {\rm Alt}_5$, the alternating group on $5$ letters
\cite{NN}.

\medskip

In case of non-perfect finite groups we prove the following
theorem.

\begin{theorem}
\label{th:fin}
Let $G$ be a finite $\Omega$-group. Then the graph
$\Delta_k^\Omega(G)$ is connected for every $k \geqslant
d_{G\Omega}(G)+1$. \label{normal2}
\end{theorem}

Note this is not true for $k = d_{G\Omega}(G)$, e.g.\ for when $G$
is abelian.

\begin{corollary} Let $G$ be a finite group and  $N \triangleleft G$
a normal subgroup.  Then the graph $\Delta_k^G(N)$ is connected
for every $k \geqslant d_G(N)+1$. \label{cor:normal2}
\end{corollary}

These results lead us to state the following conjecture.

\begin{quote}
 \textsf{\sl {\bf Relativised Finitary AC-Conjecture:} \/ Let $G$ be a finite
$\Omega$-group and $k = d_{G\Omega}(G) \geqslant 2$. Then two
tuples $U, V$ from $N_k(G,\Omega)$ are $AC_{G\Omega}$-equivalent
if and only if they are $AC_{\Omega Ab(G)}$-equivalent in the
abelianisation ${\rm Ab}(G) = G/[G,G]$, i.e., the connected
components of the graph $\Delta_k^\Omega(G)$ are precisely the
preimages of the connected components of the graph
$\Delta_k^\Omega({\rm Ab}(G))$.}
\end{quote}

Theorem~\ref{th:lifting} confirms the conjecture when $G =
\Omega$.

\section{Elementary properties of AC-transformations}
\label{se:2-2}

Let $G$ be an  $\Omega$-group. From now on for tuples $U,V \in
G^k$ we write $U \sim_G V$, or simply $U \sim V$,  if the tuples
$U, V$ are $AC_{G\Omega}$-equivalent in $G$.

\begin{lemma}
\label{le:AC-modulo} Let $G$ be an  $\Omega$-group, $N$ a normal
$\Omega$-subgroup of
$G$, and $\phi:G \rightarrow G/N$ the canonical epimorphism.
Suppose $(u_1, \ldots, u_k)$ and $(v_1, \ldots, v_k)$ are two
$k$-tuples of elements from $G$. If
$$
(u_1^\phi, \ldots, u_k^\phi) \sim_{G/N} (v_1^\phi, \ldots,
v_k^\phi)$$ then there are elements $m_1, \ldots, m_k \in N$ such
that
$$
(u_1, \ldots, u_k) \sim_{G} (v_1m_1, \ldots, v_km_k).$$
 Moreover,  one can use the same system of
elementary transformations {\rm (}after replacing conjugations by
elements $gN \in G/N$ by conjugations by elements $g \in G${\rm )}
.
\end{lemma}
{\it Proof.} Straightforward. \hfill $\Box$

\begin{lemma}
\label{le:AC-extension} Let $G$ be an $\Omega$-group. If  $(w_1,
\ldots, w_k) \in G^k$ then for every $i$ and every  element $g \in
gp_{G\Omega}(w_1, \ldots, w_{i-1}, w_{i+1}, \ldots , w_k)$
 $$(w_1, \ldots, w_k) \ \sim_G \ (w_1, \ldots, w_ig, \ldots, w_k).$$
\end{lemma}
{\it Proof.} Obvious. \hfill $\Box$

\section{The N-Frattini subgroup and semisimple decompositions}
\label{se:1}

\begin{definition}
\label{de:nf} Let\/ $G$ be an $\Omega$-group. The N-Frattini subgroup of\/
$G$
is the intersection of all proper maximal normal $\Omega$-subgroups of\/
$G$,
if such exist, and the group $G$, otherwise.  We denote it by
$W(G)$.
\end{definition}

Observe, that if $G$ has a non-trivial finite $\Omega$-quotient then $W(G)
\neq G$.

 An element $g$ in an $\Omega$-group $G$ is called {\it non-N-generating} if
for
every subset $Y \subset G$  if $gp_G(Y \cup \{g\}) = G$ then
$gp_G(Y) = G$.

\begin{lemma} \
\label{le:fr-analog}
\begin{enumerate}
\item [{\rm (1)}] The set of all non-N-generating
elements of an $\Omega$-group $G$
coincides with $W(G)$.

\item [{\rm (2)}] A tuple $U = (u_1,\ldots, u_k)$
generates $G$ as a normal $\Omega$-subgroup if and only if
the images $(\bar u_1,\ldots, \bar u_k)$ of elements $u_1,\ldots, u_k$
in $\bar G = G/W(G)$  generate $\bar G$ as normal\/ $\Omega$-subgroup.

\item [{\rm (3)}] $G/W(G)$ is an $\Omega$-subgroup of
an {\rm (}unrestricted{\rm )} Cartesian product of\/ $\Omega$-simple
$\Omega$-groups
{\rm (}that is, $\Omega$-groups which do
not have proper non-trivial normal $\Omega$-subgroups{\rm )}.

\item [{\rm (4)}] As an abstract group,  $G/W(G)$ is a subgroup of
an {\rm (}unrestricted{\rm )} Cartesian product of\/
characteristically simple groups. In particular, if $G$ is finite
then $G/W(G)$ is a product of simple groups.

\end{enumerate}
\end{lemma}
 {\it Proof.} (1) and (2) are similar to the standard proof
for the analogous property of the Frattini  subgroup.

To prove (3)  let $N_i$, $i \in I,$ be the set of all maximal proper
normal $\Omega$-subgroups of $G$. The canonical epimorphisms $G \rightarrow
G/N_i = G_i $ give rise to a homomorphism $\phi:G \rightarrow
\overline{\prod}_{i \in I}G_i$ of $G$ into the unrestricted
Cartesian product of $\Omega$-groups $G_i$. Clearly, $\ker \phi = W(G)$. So
$G/W(G)$ is an $\Omega$-subgroup of the Cartesian product of $\Omega$-simple
$\Omega$-groups $G_i$.

To prove (4) it suffices to notice that $G_i = G/N_i$ has
no $\Omega$-invariant normal subgroups, hence is characteristically
simple. \hfill  $\Box$

\medskip

To study the quotient $G/W(G)$  we need to recall
 a few definitions.
 Let
 $$G = \prod_{i \in I}G_i$$
  be a direct  product of $\Omega$-groups.   Elements $g \in G$ are
functions $g:I \rightarrow \bigcup G_i$ such that $g(i) \in G_i$ and
with finite support $supp(g) = \{i \in I \mid g_i \neq 1\}$. By
$\pi_i: G \rightarrow G_i$ we denote the canonical projection
$\pi_i(g) = g(i)$, we also denote $\pi_i(g) = g_i$. Sometimes we
identify the group $G_i$ with its image in $G$ under the canonical
embedding $\lambda_i: G_i \rightarrow G$ such that
$\pi_i(\lambda_i(g)) = g$ and $\pi_j(\lambda_i(g)) = 1$ for $j
\neq i$.

 An embedding (and we can always assume it is an
inclusion) of an $\Omega$-group $H$ into the $\Omega$-group $G$
 \begin{equation}
 \label{eq:subdirect}
 \phi : H \hookrightarrow \prod_{i \in I} G_{i}
  \end{equation}
  is called a {\it subdirect decomposition} of $H$ if
  $\pi_i(H) = G_i$ for each $i$ (here $H$
is viewed as a subgroup of $G$). The subdirect decomposition
(\ref{eq:subdirect}) is termed {\it minimal} if $H\cap G_{i}\not = \{
1 \}$ for any $i = 1, \ldots, n$, where both $G_i$ and $H$ are
viewed as subgroups of $G$. It is easy to see that given a
subdirect decomposition of $H$ one can obtain a  minimal one by
deleting non-essential factors (using Zorn's lemma).

\begin{definition}
An\/ $\Omega$-group\/ $G$ admits a finite semisimple decomposition if\/
$W(G) \neq
G$ and\/  $G/W(G)$ is a finite direct product of\/ $\Omega$-simple
$\Omega$-groups.
\end{definition}

The following lemma shows that any minimal subdirect decomposition
into simple groups is, in fact, a direct decomposition.

\begin{lemma} \label{le:dir} Let $\phi: G \rightarrow \prod_{i \in
I}G_i $ be  a minimal  subdirect decomposition of an $\Omega$-group
$G$ into $\Omega$-simple
$\Omega$-groups $G_i$, $i \in I$.   Then $G = \prod_{i \in I}G_i.$
\end{lemma}
{\it Proof.} Let $K_i = G\cap G_i$, $i \in I$.  It suffices to
show that $K_i = G_i$. Indeed, in this event $G \geqslant \prod_{i
\in I}G_i$ and hence $G = \prod_{i \in I}G_i.$.

Fix an arbitrary $i \in I$.  Since $\phi$ is minimal there exists
a non-trivial $g_i \in K_i$. For an arbitrary $x_i \in G_i$ there
exists an element $x \in G$ such that $\pi_i(x) = x_i$. It follows
that $g_i^x = g_i^{x_i} \in K_i$. Hence $K_i \geqslant
gp_{G_i\Omega}(g_i) = G_i$, as required. \hfill $\Box$

\begin{lemma}
\label{le:unique-semisimple} If an $\Omega$-group $G$ has a finite
semisimple decomposition then it is unique {\rm (}up to a permutation of
factors{\rm )}.
\end{lemma}
{\it Proof.} Obvious. \hfill $\Box$

\medskip

Obviously, an $\Omega$-group $G$ admits a finite semisimple decomposition if
and only if $W(G)$ is intersection of finitely many maximal normal
 $\Omega$-subgroups of $G$. This implies the following lemma.

 \begin{lemma}
  A finite $\Omega$-group admits a finite semisimple decomposition.
\end{lemma}

\section{Connectivity of Andrews-Curtis graphs of perfect finite
groups} \label{se:3}

Recall that a group $G$ is called perfect if $[G,G] = G$.

\begin{lemma}
Let an $\Omega$-group $G$ admits a finite semisimple decomposition:
$$G/W(G) = G_1 \times \cdots \times G_k.$$
 Then $G$ is perfect if and only if all $\Omega$-simple $\Omega$-groups
$G_i$ are
 non-abelian.
\end{lemma}
{\it Proof.} Obvious. \hfill $\Box$

\medskip

We need the following notations to study normal generating tuples
in an $\Omega$-group $G$ admitting finite semisimple decomposition. If $g
\in \prod_{i \in I} G_i$ then by $supp(g)$ we denote the set of
all indices $i$ such that $\pi_i(g) \neq 1$.
\begin{lemma}
\label{le:5} Let\/ $G = \prod_{i \in I}G_i$ be a finite product of
$\Omega$-simple non-abelian $\Omega$-groups. If\/ $g \in G$ then
$gp_{G\Omega}(g) \geqslant G_i$ for any $i \in supp(g)$.
\end{lemma}
{\it Proof.} If $g \in G$ and $g_i = \pi_i(g) \neq 1$, then there
exists $x_i \in G_i\Omega$ with $[g_i,x_i] \neq 1$. Hence $1 \neq
[g,x_i] = [g_i,x_i]  \in gp_{G\Omega}(g) \cap G_i$.  Since $G_i$
is $\Omega$-simple it coincides with the nontrivial normal
$\Omega$-subgroup $gp_{G\Omega}(g) \cap G_i$, as required. \hfill
$\Box$

\medskip

 Let $G/W(G) = \prod_{i\in I} G_i$ be the canonical semisimple
 decomposition of an $\Omega$-group $G$.
 For an element $g \in G$ by ${\bar g}$ we denote the canonical image
$gW(G)$ of
$g$ in $G/W(G)$ and by   $supp(g)$ we denote the support
$supp({\bar g})$ of ${\bar g}$. \hfill $\Box$

 \begin{lemma}
\label{co:suppdirW} Let\/ $G$ be a finite perfect $\Omega$-group
and\/ $G/W(G) = \prod_{i\in
I} G_i$ be its canonical semisimple decomposition.
Then a finite set of
elements  $g_1, \ldots, g_m \in G$ generates $G$ as
a normal $\Omega$-subgroup if and only
if
 $$supp(g_1) \cup \cdots \cup supp(g_m) = I.$$
  \end{lemma}

{\it Proof}. It follows from Lemma \ref{le:5} and Lemma \ref{le:fr-analog}.
\hfill $\Box$

\medskip

\paragraph{Proof of Theorem~\ref{normal}.} We can now prove
Theorem~\ref{normal} which settles the Relativised Finitary
AC-Conjecture in affirmative for finite perfect $\Omega$-groups.

Let $G$ be a finite perfect $\Omega$-group, $\overline{G} =
G/W(G)$,  and $\overline{G}  = \prod_{i\in I} G_i$ be its
canonical semisimple decomposition. Fix an arbitrary $k \geqslant 2$.

{\sc Claim 1.}  Let $U = (u_1, \ldots, u_k) \in N_k(G,\Omega)$. Then there
exists an element $g \in G$  with $supp(g) = I$ such that
 $$(u_1, \ldots, u_k) \sim_G (g, u_2, \ldots, u_k) .$$

Indeed, by Lemma \ref{le:fr-analog} the tuple $U$ generates $G$ as
a normal subgroup if and only if its image ${\overline U}$ generates
${\overline G}$ as a normal subgroup. Lemma \ref{le:AC-modulo}  shows
that it suffices to prove the claim for the $\Omega$-group ${\overline G}$
(recall that $supp(g) = supp({\bar g})$). So we can  assume that
$G = \prod_{i\in I} G_i$.
 Since  $U \in N_k(G,\Omega)$,  Lemma \ref{co:suppdirW} implies that
 $$supp(u_1) \cup \cdots \cup supp(u_k) = I.$$
Let $i \in I$ and $i \not \in supp(u_1)$. Then there exists an
index $j$ such that $i \in supp(u_j)$.  By Lemma \ref{le:5},
$gp_{G\Omega}(u_j) \geqslant G_i$. So there exists a non-trivial
$h \in gp_{G\Omega}(u_j)$ with $supp(h) = \{i\}$. By
Lemma~\ref{le:AC-extension}, $U \sim (u_1h, u_2, \ldots, u_k) =
U^* $ and $supp(u_1h) = supp(u_1) \cup \{i\}$.  Now the claim
follows by induction on the cardinality of $I\smallsetminus
supp(u_1)$. In fact, one can bound the number of elementary
AC-moves needed in Claim 1. Indeed, since $G_i$ is non-abelian
$\Omega$-simple there exists an element $x \in G\Omega$ such that
$u_j^x \neq u_j$. Then the element $h$ above can be taken in the
form $h = u_j^xu_j^{-1}$, and only  four moves are needed to
transform  $U$ into $U^*$. This proves the claim.

{\sc Claim 2.}  Every $k$-tuple $U_1 =  (g, u_2, \ldots, u_k)$
with $supp(g) = I$ is
AC-equivalent to a tuple $U_2 =  (g, 1, \ldots, 1)$.

 By Lemma \ref{co:suppdirW} $g$ generates $G$ as a normal $\Omega$-subgroup.
Now the
 claim follows from Lemma~\ref{le:AC-extension}.

{\sc Claim 3}. Every two $k$-tuples $U_2 =  (g, 1, \ldots, 1)$ and  $U_3 =
(h, 1,
\ldots, 1)$  from $N_k(G,\Omega)$ are AC-equivalent.

Indeed, $U_2$ is AC-equivalent to $(g, 1, \ldots, 1, g)$. By
Lemma~\ref{le:AC-extension}  the former one is AC-equivalent to
$(h, \ldots, 1,g)$, which is AC-equivalent to  $(h,1, \ldots, 1)$,
as required.

 The theorem follows from Claims 1, 2, and 3. \hfill $\Box$

 \section{Arbitrary finite groups}
 \label{se:arbitrary}

 \begin{lemma}
 \label{le:nd}
 Let
 \begin{equation}
 \label{eq:nd}
G = G_1 \times \cdots \times G_s \times A \end{equation}
  be a direct  decomposition of an $\Omega$-group $G$ into a product of
  non-abelian
$\Omega$-simple $\Omega$-groups $G_i, i = 1, \ldots s,$ and an
abelian $\Omega$-group $A$.   Then, assuming $G \ne 1$,
  $$d_{G\Omega}(G) = \max\{d_{A\Omega}(A),1\}.$$
\end{lemma}
{\it Proof.} Put $S(G) = G_1 \times \cdots \times G_s .$  Since
$A$ is a quotient of $G$  then $d_{G\Omega}(G) \geqslant
d_{A\Omega}(A)$. Therefore, $d_{G\Omega}(G) \geqslant
\max\{d_{A\Omega}(A),1\}$. On the other hand, if $g$ generates
$S(G)$ as a normal $\Omega$-subgroup (such $g$ exists by Lemma
\ref{co:suppdirW}) and $a_1, \ldots, a_{d_\Omega(A)}$ generate $A$
then we claim that the tuple of elements from $G$:
$$(ga_1, a_2,
\ldots, a_{d_\Omega(A)})$$ generates $G$ as a normal
$\Omega$-subgroup. Indeed, let $g = g_1 \cdots g_s$ with $1\neq
g_i \in G_i$. Since $G_i$ is non-abelian then $g_i$ is not
central in $G_i$ and hence there exists $h_i \in G_i$ such that
$[g_i,h_i] \neq 1$. It follows that if $h = h_1 \ldots h_s$ then
$[g,h] \neq 1$ and $supp([g,h]) = \{1, \ldots, n\}$. In
particular, $[g,h]$ belongs to $N = gp_{G\Omega}(ga_1, a_2,
\ldots, a_{d_\Omega(A)})$ and generates $S(G)$ as a normal
$\Omega$-subgroup. Therefore, $S(G) \subset N$ and hence $a_1,
\ldots, a_{d_\Omega(A)} \in N$, which implies that $G  = N$. This
shows that $d_{G\Omega}(G) = \max\{d_\Omega(A),1\}$, as required.

\paragraph{Proof of Theorem~\ref{th:fin}.}
 Let $G$ be a minimal counterexample to the statement of the
 theorem.   Then $G$ is not perfect.
$G$ is also non-abelian by Fact~\ref{fact:dg}.
 Put $t = d_{G\Omega}(G)$  and $k \geqslant t+1$. Let
 $M$ be a minimal non-trivial normal $\Omega$-subgroup of $G$.  It
 follows that $M \neq G$, and the theorem holds for the
 $\Omega$-group $\overline{ G} =
 G/M$. Obviously, $k >  d_{G\Omega}(G) \geqslant d_{\bar G\Omega}(\overline{G})
$,
 hence   the AC-graph
 $\Delta_k^\Omega (\overline{G})$ is connected.  Fix any tuple $(z_1,
 \ldots,z_t) \in N_t(G,\Omega)$. If $(y_1, \ldots,y_k)$ is an arbitrary
tuple from
 $N_k(G,\Omega)$ then the $k$-tuples $(\bar y_1, \ldots,\bar y_k)$ and
$(\bar z_1,
 \ldots,\bar z_t, 1,\ldots,1)$ are AC-equivalent in $\overline{ G}$. Hence
by Lemma
 \ref{le:AC-modulo} there
 are elements $m_1, \ldots, m_k \in M$ such that
 $$(y_1, \ldots,y_k) \sim (z_1m_1, \ldots, z_tm_t, m_{t+1},\ldots, m_k).$$
 We may assume that one of the elements $m_{t+1},\ldots, m_k$ in distinct
from $1$,
 say  $m_k \neq 1$. Indeed, if $m_{t+1} = \ldots =m_k=1$ then the elements
 $z_1m_1, \ldots, z_tm_t$ generate
 $G$ as a normal $\Omega$-subgroup, hence applying AC-transformations we
 can get any non-trivial element from $M$ in the place of $m_k$.
 Since $M$ is a minimal normal $\Omega$-subgroup of $G$ it follows that $M$
 is the $G\Omega$-normal closure of $m_k$ in $G$, in particular, every $m_i$
 is a product of conjugates of $m_k^{\pm 1}$. Applying
 AC-transformations we can get rid of all elements $m_i$, $i = 1,
 \ldots, m_t$, in the tuple above. Hence,
$$(z_1m_1, \ldots, z_tm_t,m_{t+1},\ldots, m_k) \sim (z_1, \ldots,
z_t,1,\ldots, 1, m_k).$$

But $(z_1, \ldots, z_t) \in N_t(G,\Omega)$,  hence
$$(z_1, \ldots, z_t,1,\ldots, m_k) \sim (z_1, \ldots, z_t,1,\ldots, 1).$$

We showed that any $k$-tuple $(y_1, \ldots, y_k) \in N_k(G,\Omega)$ is
AC-equivalent to the fixed tuple $(z_1, \ldots, z_t,1,\ldots,1)$. So the
AC-graph $\Delta_k^\Omega(G)$ is connected and $G$ is not a
counterexample. This proves the theorem. \hfill $\Box$

\section{Proof of Theorem~\ref{th:G=Omega}}

We denote by $\tilde g$ the image of $g\in G$ in the abelinisation
${\rm Ab}(G) = G/[G,G]$.

We systematically, and without specific references, use elementary
properties
 of Andrews-Curtis transformations,
 Lemmas ~\ref{le:AC-modulo} and \ref{le:AC-extension}.

Suppose  Theorem~\ref{th:lifting} is false. Consider   a
counterexample $G$ of minimal order for a given $k \geqslant
d_G(G)$. For a given $k$-tuple $(g_1,\dots, g_k) \in N_k(G)$  we
denote by ${\cal C}(g_1,\dots, g_k)$ the set
$$\{(h_1,\dots,h_k) \in N_k(G) \mid (\tilde{g}_1,\dots,\tilde{g}_k)
 \sim (\tilde{h}_1,\dots,
\tilde{h}_k) \ \& \ (g_1, \ldots, g_k) \not \sim (h_1,
\ldots,h_k)\}$$
 Put
  $${\cal D} = \{(g_1,\dots, g_k)
\in N_k(G) \mid {\cal C}(g_1,\dots, g_k) \neq \emptyset\}.$$ Then
the set ${\cal D}$ is not empty.  Consider the following subset of
${\cal D}$:
 $$ {\cal E} = \{(g_1,\dots, g_k) \in {\cal D} \mid  |gp_G( g_2,\dots,g_k )| \
 \hbox{ is  minimal  possible}\}.$$
 Finally, consider the
subset ${\cal F}$ of ${\cal E}$:
 $${\cal F} = \{(g_1,\dots, g_k) \in {\cal E} \mid \ |gp_G(g_1)|\
\hbox{ is   minimal   possible }\}$$
 In order to
prove the theorem it suffices to show that $G$ is abelian.

Fix an arbitrary tuple $(g_1, \ldots, g_k) \in {\cal F}$ and an
arbitrary tuple $(h_1,\ldots,h_k) \in {\cal C}(g_1,\ldots,g_k)$.
Denote $G_1 = gp_G(g_1)$ and $G_2 = gp_G( g_2,\dots,g_k )$.

The following series of claims provides various inductive
arguments which will be in use later.

Notice that the minimal choice of $g_1$ and $g_2,\dots,g_{k}$ can
be reformulated as

\bigskip \noindent \textsc{Claim 1.1} \emph{Let $f_1 \in G_1$,
$f_2,\dots,f_k\in G_2$ such that $(f_1,f_2,\dots,f_k) \in {\cal
C}(g_1, \ldots,g_k)$.  Then
$$
gp_G( f_1) = G_1 \hbox{ and } gp_G( f_2,\dots, f_k ) = G_2.
$$
}

\bigskip \noindent \textsc{Claim 1.2} \emph{Let $f_1 \in G$,
$f_2,\dots,f_k\in G_2$ such that $(f_1,f_2,\dots,f_k) \in {\cal
C}(g_1, \ldots,g_k)$.  Then
$$
 gp_G( f_2,\dots,f_k ) = G_2.
$$
}

\bigskip \noindent \textsc{Claim 1.3} \emph{ Let $M$ be a
non-trivial normal subgroup of $G$. Then
$$
(h_1,\dots,h_k) \sim (g_1m_1,\dots,g_{k-1}m_{k-1}, g_km_k)
$$
for some $m_1,\dots,m_k \in M$.
 }

\medskip\noindent
Indeed, obviously $$(h_1M,\dots,h_kM), (g_1M, \ldots,g_kM) \in
N_k(G/M).$$ Moreover, since
$$
(\tilde{g}_1,\dots,\tilde{g}_k) \sim (\tilde{h}_1,\dots,
\tilde{h}_k)
$$
there exists a sequence of AC-moves $t_1, \ldots, t_n$ (where each
$t_i$ is one of the transformations (1)--(4), with the specified
values of $w$ in the case of transformations (4))
 and elements $c_1, \ldots ,c_k
\in [G,G]$ such that
$$
(h_1,\ldots,h_k)t_1\cdots t_k = (g_1c_1,\ldots,g_kc_k)
$$
Therefore
$$
(h_1M,\ldots,h_kM)t_1\cdots t_k = (g_1c_1M,\ldots,g_kc_kM)
$$
Since  $c_iM \in [G/M,G/M]$ for every $i = 1,\ldots,k$ this shows
 that the images of the tuples $(h_1M,\ldots,h_kM)$ and
$(g_1M,\ldots,g_kM)$ are AC-equivalent in the abelianisation
$Ab(G/M)$. Now the claim follows from the fact that $|G/M| < |G|$
 and the assumption that $G$ is the minimal possible counterexample.

\medskip

 The following claim says that the set ${\cal
C}(g_1,\ldots,g_k)$ is closed under $\sim$.

\bigskip \noindent \textsc{Claim 1.4} \emph{
 If $(e_1,\ldots,e_k) \in {\cal C}(g_1, \ldots,g_k)$ and $(f_1,\ldots,f_k)
\sim (e_1, \ldots,e_k)$ then $(f_1, \ldots,f_k) \in {\cal C}(g_1,
\ldots,g_k)$
}

\medskip

Now we study the group $G$ in a series of claims.

\bigskip \noindent \textsc{Claim 2.} $G = G_1 \times G_2$.

\medskip
Indeed, it suffices to show that $G_1 \cap G_2 =1$.  Assume the
contrary, then $M = G_1 \cap G_2 \neq 1$ and by Claim  1.3
$$
(h_1,\dots,h_k) \sim (g_1m_1,\dots,g_{k-1}m_{k-1}, g_km_k)
$$
for some $m_1,\dots,m_k \in M$. By Claim 1.4
$$
(g_1m_1,\dots,g_{k-1}m_{k-1}, g_km_k) \in {\cal C}(g_1,\ldots,g_k)
$$
 By Claim 1.1,
$$
gp_G( g_1m_1 ) = G_1, \quad gp_G( g_2m_2,\dots, g_k m_k ) = G_2
$$
and we can represent the elements $m_2,\dots,m_k\in G_1\cap G_2$
as products of conjugates of $g_1m_1$, therefore deducing that
$$
(g_1m_1,g_2m_2\dots, g_km_k) \sim (g_1m_1,g_2,\dots, g_k).
$$
Since $m_1 \in gp_G( g_2,\dots,g_k\rangle$, we conclude that
$$
(g_1m_1,g_2, \dots,
 g_k) \sim
(g_1,g_2,\dots, g_k),
$$
and therefore
$$
(h_1,\ldots, h_k) \sim (g_1,\ldots, g_k),
$$
a contradiction. This proves the claim.
 \hfill $\square$

\bigskip \noindent \textsc{Claim 3.}
$ [G_2, G_2] = 1.$ \emph{In particular, $G_2 \leqslant Z(G)$.}

\medskip

Indeed, assume the contrary. Then $M = [G_2, G_2] \neq 1$ and  by
Claim  1.3
$$
(h_1,\dots,h_k) \sim (g_1m_1,\dots, g_km_k), \quad m_1,\dots, m_k
\in M \leqslant G_2 .
$$
By virtue of Claims 1.4 and 1.2,
 $gp_G(g_2m_2,\dots,g_km_k ) = G_2$ and hence $m_1 \in gp_G(
g_2m_2,\dots,g_km_k ).$ It follows that
$$
(g_1m_1,g_2m_2,\dots, g_km_k) \sim (g_1,g_2m_2\dots, g_km_k).
$$
Therefore it will be enough to prove
$$
(g_1,g_2m_2,\dots, g_km_k) \sim (g_1,g_2,\dots, g_k).
$$
We proceed as follows, systematically using the fact that
$g_2,\dots,g_{k}$ and all their conjugates commute with all the
conjugates of $g_1$.

We start with a series of Nielsen moves which lead to
\begin{eqnarray*}
(g_1,g_2m_2\dots, g_km_k) & \sim & (g_1,\, g_1\cdot g_2m_2,g_3m_3,\dots, g_km_k)\\
& \sim & (g_1\cdot m_2,\, g_1g_2m_2,g_3m_3,\dots, g_km_k).
\end{eqnarray*}
The last transformation is the key for the whole proof and
requires some explanation. Since $m_2$ belongs to
$$
[G_2,G_2] = [gp_G( g_2m_2,\dots, g_km_k ), gp_G( g_2m_2,\dots,
g_km_k )],
$$
$m_2$ can be expressed as a word
$$w(x_2,\dots,x_k)=(x_{i_1}^{f_1})^{\varepsilon_1} \cdots
(x_{i_l}^{f_l})^{\varepsilon_l}
$$
where  $x_i= g_im_i$, $i=2,\dots,k$,  $f_j \in G$ and the word
$w$ is balanced for each variable $x_i$, that is, for each $h =
2,\dots, k$, the sum of exponents for each $x_h$ is zero:
$$\sum_{{i_j} =h}\varepsilon_j=0.$$ Moreover, since $G = G_1\times
G_2$, we can choose $f_j \in G_2$, whence commuting with $g_1\in
G_1$. Therefore
$$
w(g_1x_2,x_3,\dots,x_k) = w(x_2,x_3,\dots,x_k)
$$
and
 $$
 w(g_1g_2m_2, g_3m_3,\dots,g_km_k) = m_2.
 $$
 Hence, by
several consecutive multiplications by appropriate conjugates of
$g_1g_2m_2$ and $g_im_i$, $i=3,\dots,k$, we can produce the factor
$m_2$ in the leftmost position in the tuple. We now continue:
\begin{eqnarray*}
(g_1m_2,\, g_1g_2m_2,g_3m_3,\dots, g_km_k) & \sim &
(g_1m_2,\,g_1g_2m_2 \cdot (g_1m_2)^{-1},g_3m_3, \dots,
g_km_k) \\
&=& (g_1m_2,\,g_2, g_3m_3,\dots, g_km_k).
\end{eqnarray*}
Again by Claims 1.4 and 1.2  $G_2 = gp_G( g_2, g_3m_3,\dots,
g_km_k\rangle$.  Since $m_2 \in G_2$,
\begin{eqnarray*}
(g_1m_2,\,g_2, g_3m_3,\dots, g_km_k) & \sim & (g_1,\,g_2,\,
g_3m_3,\dots, g_km_k).
\end{eqnarray*}
Next we want to kill $m_3$. Present $m_3$ as a balanced word in
$g_2, g_3m_3,\dots, g_km_k$ conjugated by elements $f_i \in G_2$.
Note that they all commute with $g_1$. As before,
$$
m_3 = w(g_2,g_1g_3m_3,g_4m_4,\dots, g_km_k)
$$
(and, actually, $m_3 = w(g_2,y_3,\dots, y_k)$ where $y_i$ are
arbitrarily chosen from  $g_im_i$ or $g_1g_im_i$, $i=3,\dots,k$.).

 Thus we have:
\begin{eqnarray*}
(g_1,\,g_2, g_3m_3,\dots, g_km_k) & \sim &
(g_1,g_2,g_1g_3m_3,g_4m_4,\dots, g_1g_km_k)\\
&\sim & (g_1m_3,g_2,g_1g_3m_3,g_4m_4,\dots, g_km_k)\\
&\sim & (g_1m_3,g_2,g_3,g_4m_4,\dots, g_km_k)\\
&\sim & (g_1,g_2,g_3,g_4m_4,\dots, g_km_k)
\end{eqnarray*}
(the last transformation uses the fact that
$gp_G(g_2,g_3,g_4m_4,\dots, g_km_k)=G_2$  by Claims 1.4 and 1.2).

 One can easily observe that we can continue this
argument in a similar way until we come to $(g_1,g_2,\dots,g_k)$ -
contradiction, which completes the proof of the claim. \hfill
$\square$

\bigskip \noindent \textsc{Claim 4.}
$$
[G_1, G_1] = 1. $$

Let $[G_1, G_1] \neq 1$. For a proof, take a minimal non-trivial
normal subgroup $M$ of $G$ which lies in $[G_1, G_1]$. Again, by
Claim 1.3, we conclude that
$$
(h_1,\dots,h_k) \sim (g_1m_1,g_2m_2,\dots, g_km_k)
$$
for some $m_1,\dots,m_k \in M$. We assume first that $M \leqslant
gp_G( g_1m_1)$. Then
\begin{eqnarray*}
 (g_1m_1,g_2m_2\dots, g_km_k)& \sim &  (g_1m_1,g_2,\dots,
g_k)
\end{eqnarray*}
and $gp_G( g_1m_1) = gp_G( g_1)$ by Claims 1.4 and  1.1. In
particular,
$$
M \leqslant [gp_G( g_1m_1), gp_G( g_1m_1)] = [gp_G( g_2g_1m_1),
gp_G( g_2g_1m_1)],
$$
where the last equality follows from the observation that $g_2 \in
Z(G)$.
 We shall
use this in further transformations:
\begin{eqnarray*}
(g_1m_1, g_2,\dots,g_k) & \sim & (g_1m_1, g_2g_1m_1,g_3,\dots, g_k )\\
 & \sim & (g_1, g_2g_1m_1,g_3,\dots,g_k )\\
  & \sim & (g_1, g_2,g_3,\dots,g_k).
\end{eqnarray*}
This shows that $(h_1,\ldots,h_k) \sim (g_1,\dots,g_k)$ -
contradiction.  Therefore we can assume that $M \not\subseteq
gp_G( g_1m_1)$ and hence $M \cap gp_G( g_1m_1) =1$. We claim that
not all of the elements $m_2, \dots,m_k$, are trivial. Otherwise
$$
(h_1,\dots,h_k) \sim (g_1m_1,g_2,\dots,g_k),
$$
and we can repeat the previous argument and come to a
contradiction. So we assume, with out loss of generality, that
$m_2 \ne 1$.

If $M$ is non-abelian then
$$
M = [M,M] = [gp_G( m_2 ),gp_G( m_2 )] = [gp_G( g_2m_2 ),gp_G(
g_2m_2 )]
$$
and
\begin{eqnarray*}
 (g_1m_1, g_2m_2,g_3m_3,\dots,g_km_k)& \sim &  (g_1,
g_2m_2,g_3m_3,\dots,g_km_k)\\
& \sim &  (g_1, g_2,g_3,\dots,g_k);
\end{eqnarray*}
we use in the last transformation that $gp_G(g_1) = G_1 \geqslant
M$.

 Therefore we can assume that $M$ is abelian. Since $M \cap
gp_G( g_1m_1) =1$ we conclude that $[M,gp_G( g_1m_1)] =1$. But
then $[M,gp_G( g_1)] =1$. In particular, $M \leqslant Z(G)$ and
the subgroup $[gp_G( g_1m_1), gp_G( g_1m_1)] = [gp_G( g_1), gp_G(
g_1)]$ contains $M$. But this is a contradiction with $M \cap
gp_G( g_1m_1) =1$. This proves the claim. \hfill $\square$

\bigskip\noindent\textsc{Final contradiction.} Claims 3 and 4 now
yield that $G$ is abelian, as required.
 \hfill $\square$ $\square$

\subsection*{Final comments}


The referee has kindly called to our attention that the result of
Myasnikov \cite{My} (mentioned in the Introduction) was also
proved independently in 1978 by Wes Browning (unpublished).

\bigskip

\bigskip

 \noindent \textsf{Alexandre V. Borovik, School of
Mathematics, PO Box 88, The University of Manchester, Sackville
Street, Manchester M60 1QD, United Kingdom}; {\tt
borovik@manchester.ac.uk}; {\tt http://www.ma.umist.ac.uk/avb/}

\

 \noindent \textsf{Alexander Lubotzky, Department of Mathematics, Hebrew
University,
 Givat Ram, Jerusalem 91904, Israel};
 {\tt alexlub@math.huji.ac.il}

\

\noindent \textsf{Alexei G. Myasnikov,
 Department of Mathematics, The City  College  of New York, New York,
NY 10031, USA}; {\tt alexeim@att.net}; {\tt
http://home.att.net/\~\,alexeim/index.htm}
\end{document}